\documentclass[12pt]{article}
\usepackage{amssymb}

\usepackage{amsmath}
\usepackage{amscd}
\usepackage{graphicx}

\newtheorem{theorem}{Theorem}[section]
\newtheorem{proposition}[theorem]{Proposition}
\newtheorem{lemma}[theorem]{Lemma}
\newtheorem{corollary}[theorem]{Corollary}
\newtheorem{conclusion}[theorem]{Conclusion}

\begin{document}

\title{Infinitesimal deformation quantization of complex analytic spaces}
\author{V. P. Palamodov\footnote{Tel Aviv University, Ramat Aviv 69978, Israel; e-mail: palamodo@tau.ac.il}}
\date{31.01.06}
\maketitle

\section{Introduction}

For the physical aspects of deformation quantization we refer to the
expository and survey papers \cite{DS} and \cite{S}. Our objective is to
initiate a version of global theory of quantization deformation in the
category of complex analytic spaces in the same lines as the theory of
(commutative) deformation. The goal of inifinitesimal theory is to do few
steps towards construction of a star-product in the structure sheaf $%
\mathcal{O}_{X}$ of holomorphic functions on a complex analytic space $X$,
occasionally with singularity. A formal power series with representing a
star-product can have a chance to converge on holomorphic functions for, at
least, a sequence of values of the parameter. This phenomenon is described
in Berezin's global theory, \cite{Be}, \cite{Be2}. We shall see strong
similarity between commutative and skew-commutative deformations, in
particular, any Poisson bracket and any Kodaira-Spencer class (commutative
deformation) are locally in ''one flacon'', that is in the same Hochschild
cohomology space. On the global level, extension of an arbitrary
infinitesimal quantization meets obstructions which are elements of the \v{C}%
ech cohomology of the sheaves of the \textit{analytic }Hochschild cohomology
in the same way as infinitesimal deformations of complex analytic spaces do.
The construction of analytic Hochschild (co)homology is given in terms of
analytic tensor products and bounded linear mappings of analytic algebras.
We prove here the general result: the analytic Hochschild (co)homology on a
complex analytic space is always a coherent analytic sheaf in each degree.

For K3-surfaces a global infinitesimal quantization is explicitly written in
terms of Jacobian brackets.

\section{Quantization of an algebra}

Let $\mathbf{k}$ be a fields of zero characteristic, $A$ be an associative
commutative $\mathbf{k}$-algebra. A \textit{quantization} of $A$ is a
bilinear operation in $A$ of the form 
\begin{equation}
a\ast b=ab+\lambda p\left( a,b\right) +\lambda ^{2}p_{2}\left( a,b\right)
+\lambda ^{3}p_{3}\left( a,b\right) +...  \label{2}
\end{equation}
where $\lambda $ is a ''small'' parameter taking values in $\mathbf{k}$ and $%
p,p_{2},p_{3},...$ are some bilinear mappings $A\times A\rightarrow A$. The
mapping $p$ is skew-symmetric and is choisen in advance (Poisson bracket of
a deterministic system). This operation, called star-product, must be an
associative operation, that is $\left( a\ast b\right) \ast c=a\ast \left(
b\ast c\right) $ for any $a,b,c\in A$. Substituting (\ref{2}) yields 
\begin{align}
& \left( a\ast b\right) \ast c-a\ast \left( b\ast c\right) =\lambda \left[
-\partial p\left( a,b,c\right) \right]  \notag \\
& +\lambda ^{2}\left[ p\left( p\left( a,b\right) ,c\right) -p\left(
a,p\left( b,c\right) \right) -\partial p_{2}\left( a,b,c\right) \right]
+\lambda ^{3}\left[ \cdot \right] +...  \label{30}
\end{align}
The equation for arbitrary $a,b,c\in A$%
\begin{equation}
\partial p\left( a,b,c\right) \doteq ap\left( b,c\right) -p\left(
ab,c\right) +p\left( a,bc\right) -p\left( a,b\right) c=0  \label{1}
\end{equation}
must be fulfilled, which is equivalent to vanishing of the linear term in (%
\ref{30}). Consider the Hochschild cochain complex of the algebra $A:\mathrm{%
C}^{\ast }\left( A\right) =\oplus _{n\geq 0}\mathrm{C}^{n}\left( A\right) $,
where $\mathrm{C}^{n}\left( A\right) $ is the space of $\mathbf{k}$%
-polylinear operators $h:A^{\otimes n}\rightarrow A$ (or of $A$-morphisms $%
A\otimes A^{\otimes n}\rightarrow A$) with the standard ''bar''-
differential $\partial .$ The cohomology of this complex is called
Hochschild cohomology and denoted $\mathrm{Hoch}^{\ast }\left( A,A\right) .$
The condition (\ref{1}) means that $p$ is a 2-cocycle. Any cocycle defines
an element $\mathrm{cl}\left( p\right) \in \mathrm{Hoch}^{2}\left(
A,A\right) $ and any cohomology class contains only one skew-symmetric
cocycle, since any coboundary element $\partial g\left( a,b\right) =ag\left(
b\right) -g\left( ab\right) +g\left( a\right) b\;$is symmetric. Denote by $%
\mathrm{Q}\left( A\right) $ the submodule of $\mathrm{Hoch}^{2}\left(
A,A\right) $ of all skew-symmetric cocycles $p$.

The Harrison cohomology space $\mathrm{Harr}^{n}\left( A,A\right) $ is the $%
\mathbf{k}$-subspace of $\mathrm{Hoch}^{n}\left( A,A\right) $ generated by
operators $h:A^{\otimes n}\rightarrow A$ that vanish on all shuffle\textit{\ 
}products of chains, see \cite{Ba}. In particular, the second order \textit{%
Harrison } cohomology $\mathrm{Harr}^{2}\left( A,A\right) $ is the subspace
(submodule) in $\mathrm{Hoch}^{2}\left( A,A\right) $ generated by symmetric
cocycles $s:A\otimes A\rightarrow A.$ We have 
\begin{equation*}
\mathrm{Hoch}^{2}\left( A,A\right) =\mathrm{Harr}^{2}\left( A,A\right)
\oplus \mathrm{Q}\left( A\right) ,
\end{equation*}
since any 2-cocycle is the sum of its symmetric and skew-symmetric parts.
The first term is the group of first infinitesimal commutative deformations
and the second one is that for deformation quantizations.

The third cohomology groups are responsible for obstructions:

\begin{proposition}
\label{obst}Suppose that $\mathrm{Harr}^{3}\left( A,A\right) =0.$ Then the
Jacobi equation 
\begin{equation*}
p\left( p\left( a,b\right) ,c\right) +p\left( p\left( b,c\right) ,a\right)
+p\left( p\left( c,a\right) ,b\right) =0
\end{equation*}
is sufficient and necessary for existence of a mapping $p_{2}:A\otimes
A\rightarrow A$ such that the second term in (\ref{30}) vanishes. In other
words, the Jacobi sum is the first obstruction to extension of the cocycle
to a star-product.
\end{proposition}

\textsc{Proof.} Take the elements in the group algebra $\mathbb{Q}\left(
S_{3}\right) $ 
\begin{align*}
e_{1}\left( a,b,c\right) & =\frac{1}{6}\left[ 2\left( a,b,c\right) -2\left(
c,b,a\right) +\left( a,c,b\right) -\left( b,c,a\right) +\left( b,a,c\right)
-\left( c,a,b\right) \right] , \\
e_{2}\left( a,b,c\right) & =\frac{1}{2}\left[ \left( a,b,c\right) +\left(
c,b,a\right) \right] , \\
e_{3}\left( a,b,c\right) & =\frac{1}{6}\left[ \left( a,b,c\right) -\left(
c,b,a\right) -\left( a,c,b\right) +\left( b,c,a\right) -\left( b,a,c\right)
+\left( c,a,b\right) \right] , \\
e_{1}+e_{2}+e_{3}& =1.
\end{align*}
They are projectors in the Hochschild complex and commute with the
differential. According to \cite{Ba} and \cite{GS}, $\mathrm{Hoch}^{3}=%
\mathrm{Hoch}_{1}^{3}\oplus \mathrm{Hoch}_{2}^{3}\oplus \mathrm{Hoch}%
_{3}^{3},$ where$\;\mathrm{Hoch}_{1}^{3}\,\cong \mathrm{Harr}^{3}.$ The
element $q\left( a,b,c\right) \doteq p\left( p\left( a,b\right) ,c\right)
-p\left( a,p\left( b,c\right) \right) $ is a $3$-cocycle and we have $%
q=e_{1}q+e_{2}q+e_{3}q,$ where all the terms are cocycles and the term $%
e_{1}q$ is a cocycle in the Harrison cohomology. It is a coboundary since of
the assumption. Moreover, 
\begin{equation*}
2e_{2}q\left( a,b,c\right) =p\left( p\left( a,b\right) ,c\right) -p\left(
a,p\left( b,c\right) \right) +p\left( p\left( c,b\right) ,a\right) -p\left(
c,p\left( b,a\right) \right) =0,
\end{equation*}
since $p$ is skew-symmetric. The last term is equal to the sum 
\begin{equation*}
e_{3}q\left( a,b,c\right) =\frac{2}{3}\left[ p\left( p\left( a,b\right)
,c\right) +p\left( p\left( b,c\right) ,a\right) +p\left( p\left( c,a\right)
,b\right) \right]
\end{equation*}
where the Jacobi bracket appears. $\blacktriangleright $

\textbf{Definition. }A skew-symmetric 2-cocycle $p$ that satisfies the
Jacobi equation is called \textit{Poisson }bracket in $A$.

\section{Cohomology of analytic algebras}

For a natural $n$ the an algebra $R_{n}=\,\mathbb{C}\left\{
z_{1},...,z_{n}\right\} \;$of convergent power series of $n$ variables is
called regular analytic $\mathbb{C}$-algebra. An \textit{analytic}
algebra is by definition a $\mathbb{C}$-algebra $A$ that admits a
isomorphism of $\mathbb{C}$-algebra\ $A\cong R_{n}/I$ \ for some $n$ and an
ideal $I$ in $R_{n}$.  Let $f_{1},...,f_{m}$ be a set of generators of $I$
(an arbitrary ideal in $R_{n}$ has finite set of generators). The algebra $A$
is the structure algebra of the germ of the analytic space 
\begin{equation*}
X=\left\{ z\in \mathbb{C}^{n},f_{1}\left( z\right) =...=f_{m}\left( z\right)
=0\right\} 
\end{equation*}
at a marked point $\bullet \in X;$ in other notations $A=\mathcal{O}\left(
X,\bullet \right) .$ 

We adapt the cohomology theory for analytic algebras taking the analytic
Hochschild homology and cohomology instead of the algebraic ones. For this
we introduce bornology in an analytic algebra and use construction of
analytic tensor product. A bornological $\mathbb{C}$-space is a vector space 
$U$ with a family $\left\{ B\right\} $ of subsets called \textit{%
bounded} sets. The family $\left\{ B\right\} $ is subjected to several natural
conditions, which hold for the bornology of arbitrary locally convex topology. 
For arbitrary bornological vectors spaces $U,V$ a linear
operator $\alpha :U\rightarrow V$ is called \textit{bounded}, if the image
of an arbitrary bounded set $B\subset U$ is bounded in $V.$ Any analytic
algebra $A=R_{n}/I$ possesses the canonical linear bornology: a set $%
B\subset R_{n}$ is bounded, if all the series $a\in B$ converge in the ball $%
\left\{ z;\left| z\right| <\varepsilon \right\} $ for some $\varepsilon >0$
and the sums are bounded by the same constant $C$ in this ball. A set $%
B^{\prime }$ is bounded in $A,$ if it is contained in the image of a bounded
set  $B\subset R_{n}.$ Any ideal $I$ is sequentially in this bornology and
the bornology of $A$ is separated. The \textit{analytic} tensor product of
regular algebras is, by definition, $R_{m}\odot R_{n}=R_{n+m}.$ For
arbitrary analytic algebras $A=R_{m}/I,B=R_{n}/J$ the analytic tensor
product is the analytic algebra $A\odot B=R_{m+n}/\left( I,J\right) ,$ where 
$\left( I,J\right) $ is the ideal generated by images of elements of $I$ and 
$J$ in $R_{m+n}.$ For an analytic algebra $A$ the chain complex $\mathrm{Ca}%
_{\ast }\left( A\right) $ is the direct sum of analytic algebras $\mathrm{Ca}%
_{k}\left( A\right) \doteq A^{\odot k+1}.$ The chain differential $\partial $
is well defined in the module $\mathrm{Ca}_{\ast }\left( A\right) ;$ it is
the alternate sum of diagonal mappings $\partial _{i}:A^{\odot
k+1}\rightarrow A^{\odot k},i=0,...,k.$ On the level of a regular algebra $%
R_{n}\rightarrow A$ the operator $\partial _{i},i=0,...,k-1$ is the diagonal
mapping: $c\left( z^{0},z^{1},...,z^{k}\right) \mapsto c\left(
z^{0},...,z^{i},z^{i},...z^{k-1}\right) ,$ where $z^{0},...,z^{k}\in \mathbb{%
C}^{n},$ and $\partial _{k}:c\mapsto c\left( z^{0},...,z^{k-1},z^{0}\right) .
$ These operators are bounded with respect to the bornology in the algebras $%
A^{\odot k+1}$. The homology of the complex $\left( \mathrm{Ca}_{\ast
}\left( A\right) ,\partial \right) $ is called \textit{analytic} Hochschild
homology; we keep the notation $\mathrm{Hoch}_{\ast }\left( A,A\right) $.
The analytic chain complex and its homology are $A$-modules: the action of
an element $a\in A$ is defined by $a\left( z\right) \cdot c\left(
z^{0},...,z^{k}\right) =a\left( z^{0}\right) c\left( z^{0},...,z^{k}\right) .
$

The cochain space $\mathrm{Ca}^{k}\left( A\right) $ is defined to be the
module of all bounded $A$-morphisms $h:A^{\odot k+1}\rightarrow A$
with respect to the standard linear bornology of analytic algebras. We call
the cohomology of the complex $\left( \mathrm{Ca}^{\ast }\left( A\right)
,\partial ^{\ast }\right) $ the \textit{analytic} Hochschild cohomology, for
which we use the same notation $\mathrm{Hoch}^{\ast }\left( A,A\right) .$

From now on we shall use the same notation $\mathrm{Hoch}^{\ast }$, $\mathrm{%
Harr}^{\ast }$ and $\mathrm{Q}$ for cohomology of the analytic cochain
complex $\mathrm{Ca}^{\ast }\left( A\right) $.

\textbf{Remark.} Any differential operator $D:A^{\otimes k}\rightarrow A$ is
a linear bounded mapping and can be extended to a bounded $A$-morphism $%
A^{\odot k+1}$ by $D\left( a\left[ \cdot \right] \right) =aD\left( \left[
\cdot \right] \right) ,$ $a\in A.$ The inverse is not true: there are
bounded mappings which are not differential operators, see Lemma \ref{s} below. 

\begin{proposition}
\label{Q}For any formal or analytic algebra $A=R/\left(
f_{1},...,f_{m}\right) $ as above the $\mathbf{k}$-space of bounded
skew-symmetric cocycles is isomorphic to 
\begin{equation*}
\mathrm{Q}\left( A\right) \cong \mathrm{Ker}\left\{ J:A^{n}\wedge
A^{n}\rightarrow A^{nm}\right\} ,
\end{equation*}
where $J$ is the mapping generated by the Jacobian $m\times n$-matrix $%
J=\left\{ \partial _{i}f_{j}\right\} .$
\end{proposition}

\begin{lemma}
\label{d}Any cocycle $p\in \mathrm{Q}\left( A\right) $ is a derivation in
each arguments.
\end{lemma}

\textsc{Proof}\textrm{.} Take an arbitrary skew-symmetric operator $p$ that
fulfils (\ref{1}) and check is a derivation with respect to each argument.
Applying (\ref{1}) to the arguments $\left( b,c,a\right) $ and to $\left(
b,a,c\right) $ yields 
\begin{equation*}
ap\left( b,c\right) =bp\left( c,a\right) -p\left( bc,a\right) +p\left(
b,ca\right) ,\,p\left( ab,c\right) =bp\left( a,c\right) +p\left( b,ac\right)
-p\left( b,a\right) c.
\end{equation*}
Substituting in (\ref{1}) yields $0=2\left[ bp\left( c,a\right) -p\left(
bc,a\right) +cp\left( b,a\right) \right] $ that is $p\left( bc,a\right)
=bp\left( c,a\right) +cp\left( b,a\right) $ that is $p$ is a derivation in
the first argument. It is also a derivation in the second argument since of
symmetry. Lower $p$ to a skew-symmetric mapping $q:R\otimes R\rightarrow A.$
It is again a derivation and is \textit{bounded} in both arguments.
Therefore it can written in the form 
\begin{equation*}
q_{R}\left( a,b\right) =\sum_{i<j}q^{ij}\left( \partial _{i}a\partial
_{j}b-\partial _{j}a\partial _{i}b\right) \,
\end{equation*}
for some $q^{ij}\in A,$ (which is easy to check by means of Leibniz
formula). Vice versa any such bivector field can be lifted to a mapping $p,$
if and only if $q_{R}\left( fa,b\right) \in I$ for any $a,b\in R,f\in I.$
This is equivalent to the condition $q_{R}\left( f,b\right) \in I$ for any $%
b\in A$ and $f\in I.$ Taking $b=z_{i}$ yields 
\begin{equation*}
\sum_{j}q^{ij}\partial _{j}f\in I,i=1,...,n,\,f\in I
\end{equation*}
that is $J\left\{ q^{ij}\right\} \in I^{nm}$ and vice versa. $%
\blacktriangleright $

\begin{corollary}
Let $\left( X,\mathcal{O}_{X}\right) $ be an arbitrary complex analytic
space. The sheaf $\mathcal{Q}_{X}$ of spaces $\mathrm{Q}\left( \mathcal{O}%
_{X,x}\right) \;$is a coherent analytic sheaf.
\end{corollary}

Compare this result with parallel fact on commutative deformations. For any
analytic algebra $A$ we have a natural isomorphism $\mathrm{Harr}^{\ast
+1}\left( A,A\right) \cong T^{\ast }\left( A\right) ,$ where $T^{\ast }$ is
the Tyurina cohomology of analytic algebra, see \cite{Pa}. The ideal $%
I\subset R_{n}$ and the algebra $A=R_{n}/I$ are called \textit{complete
intersection, }if \textrm{dim\thinspace }$X=n-m$ for a set of generators $%
f_{1},...,f_{m}$ of $I.$

\begin{proposition}
The equations $T^{i}\left(
A\right) =0,i>1$ hold for an arbitrary complete intersection analytic algebra $A$.
\end{proposition}

This implies that \textrm{Harr}$^{3}\left( A,A\right) =0$ for any complete
intersection algebra $A$ and by Proposition \ref{obst} the mapping $p_{2}$
can be found for any Poisson bracket in $A$.

\section{Jacobians as Poisson brackets}

Take $A=R_{3}/\left( f\right) ,\,X=\left\{ z:f\left( z\right) =0\right\} .$
The set $X$ is a complex analytic surface, occasionally singular one. The
module $\mathrm{Q}\left( A\right) $ contains the Jacobian bracket defined on 
$R_{3}$ by 
\begin{equation*}
P_{f}\left( a,b\right) =J_{x}\left( f,a,b\right) =\mathrm{\det }\left( 
\begin{array}{ccc}
\partial _{1}f & \partial _{2}f & \partial _{3}f \\ 
\partial _{1}a & \partial _{2}a & \partial _{3}a \\ 
\partial _{1}b & \partial _{2}b & \partial _{3}b
\end{array}
\right) ,
\end{equation*}
where $\partial _{i}=\partial /\partial z_{i}.$ It can be lifted to the
bracket $P_{f}:A\otimes A\rightarrow A$, since 
\begin{equation*}
P_{f}\left( fa,b\right) =J_{x}\left( f,fa,b\right) =fJ_{x}\left(
f,a,b\right) +aJ_{x}\left( f,f,b\right) .
\end{equation*}

\begin{proposition}
\label{jac}The mapping $eP$ fulfils the Jacobi identity for any $e\in A$.
\end{proposition}

\textrm{Proof. }The factor $e$ can be eliminated replacing $f$ by $ef.$
First we check the identity for arbitrary $a,b,c,d\in A$ 
\begin{equation*}
\sum_{a\mapsto b\mapsto c\mapsto a}J\left( f,a,b\right) J\left( f,c,d\right)
=0
\end{equation*}
where the sum is taken over cyclic permutations of $a,b,c$. We can assume
that \textrm{d}$f\neq 0$ since the equation is purely algebraic. Then we
change the variables to $y_{1},y_{2},y_{3}=f\left( x\right) $ and cancel out
the factor \textrm{det}$^{2}\partial y/\partial x.$ Then we have $%
J_{y}\left( f,a,b\right) =\partial _{1}a\partial _{2}b-\partial
_{2}a\partial _{1}b$ and the proof is straightforward. The is true for the
identity 
\begin{equation*}
\sum_{a\mapsto b\mapsto c\mapsto a}J_{y}\left( f,a,J_{y}\left( f,b,c\right)
\right) =0
\end{equation*}
We have 
\begin{eqnarray*}
J_{x}\left( f,a,J_{x}\left( f,b,c\right) \right) &=&dJ_{y}\left(
f,a,dJ_{y}\left( f,b,c\right) \right) \\
&=&d^{2}\left[ J_{y}\left( f,a,J_{y}\left( f,b,c\right) \right) +J\left(
f,a,d\right) J\left( f,b,c\right) \right]
\end{eqnarray*}
where $d\doteq \;$det$\partial y/\partial x.$ The sum over cyclic
permutation of the right-hand side vanishes, since of the above arguments.$%
\blacktriangleright $

\textbf{Example. }Take the function $f\left( x\right) =1/2\left(
x_{1}^{2}+x_{2}^{2}\pm x_{3}^{2}\right) $ in $\mathbb{R}^{3}$ and consider
the Poisson bracket $J\left( f,a,b\right) $. It defines the Poisson
structure in the algebra of motions of solid 3-body with a fixed point or in
the real quadratic cone, depending on the sign.

Let $A$ $=R_{n}/(f_1,...,f_m)$ be a complete intersection analytic algebra. 
Consider the Jacobian matrix 
\begin{equation*}
J\left( f_{1},...,f_{m},a,b\right) \doteq \det \left( 
\begin{array}{cccc}
\partial _{1}f_{1} & \partial _{2}f_{1} & ... & \partial _{n}f_{1} \\ 
... & ... & ... & ... \\ 
\partial _{1}f_{m} & \partial _{2}f_{m} & ... & \partial _{n}f_{m} \\ 
\partial _{1}a & \partial _{2}a & ... & \partial _{n}a \\ 
\partial _{1}b & \partial _{2}b & ... & \partial _{n}b
\end{array}
\right) ,a,b\in R_{n}
\end{equation*}
Let $K$ be a subset of $\left[ 1,...,n\right] $ of $m+2$ elements and $%
J_{K}\left( f_{1},...,f_{m},a,b\right) $ be the corresponding minor of this
matrix.

\begin{proposition}
\label{pb}For any set $K$ as above the mapping 
\begin{equation*}
P_{K}\left( a,b\right) =J_{K}\left( f_{1},...,f_{m},a,b\right) 
\end{equation*}
can be lifted to a bilinear mapping $p_{K}:A\otimes A\rightarrow A$. This is
a Poisson bracket. Moreover, any $\mathbf{k}$-linear combination of brackets 
$p_{K}$ is a Poisson bracket. \newline
If the germ  $X\backslash \cdot $ is regular any skew-symmetric 2-cocycle $q$ is
generated by such brackets: 
\begin{equation*}
q=\sum_{K\subset \left[ 1,...,n\right] }e_{K}P_{K},\;e_{K}\in A.
\end{equation*}
\end{proposition}

\textbf{Remark. }The sum $q$ need not to satisfy the Jacobi identity,
whereas each term does.

\section{Globalization and obstructions}

Let $\left( X,\mathcal{O}_{X}\right) $ be a complex analytic space and $%
\mathcal{Q}_{X}$ be the sheaf of germs of skew-symmetric 2-cocycles on $X.$
Any section $p\in \Gamma \left( X,\mathcal{Q}_{X}\right) $ can be considered
as first infinitesimal of a global quantization of $X$ that is a first order
term in a global star-product of the form (\ref{2}), where $a,b$ are
arbitrary elements of the structure sheaf $\mathcal{O}_{X}$ and $p_{k}:%
\mathcal{O}_{X}\times \mathcal{O}_{X}\rightarrow \mathcal{O}_{X},\,k=2,3,...$
are some bilinear bounded mappings. The term $p_{2}$ has to fulfil the
cohomological equation 
\begin{equation}
\left[ p,p\right] \left( a,b,c\right) \doteq p\left( p\left( a,b\right)
,c\right) -p\left( a,p\left( b,c\right) \right) =-\partial p_{2}\left(
a,b,c\right)   \label{5}
\end{equation}
There are some obstructions that can be revealed in several steps:

(i) The cohomology class of the Gerstenhaber bracket $\left[ p,p\right] $ is
an obstruction to existence of $p_{2}.$ The arguments of Sec. 3 show that
the equation $e_{3}\left[ p,p\right] =0$ is satisfied if and only if $p$ is
a Poisson bracket and the equation $e_{2}\left[ p,p\right] =0$ is always
fulfilled. The cohomology class $e_{1}\left[ p,p\right] $ belongs for any
point $x\in X$ to the Harrison cohomology $\mathrm{Harr}^{3}\left( \mathcal{O%
}_{X,x},\mathcal{O}_{X,x}\right) .$ In this step, we assume that the mapping 
$p$ is locally \textit{bounded} in the sense of Sec. 4. The bracket and the
cocycle $e_{1}\left[ p,p\right] $ are also bounded. Therefore the class $%
\mathrm{cl}\left( e_{1}\left[ p,p\right] \right) $ is contained in the
analytic Harrison cohomology, which is isomorphic to Tyurina cohomology $%
T^{2}\left( \mathcal{O}_{X,x}\right) .$ This module is the stalk of the
coherent sheaf $\mathcal{T}^{2}\left( \mathcal{O}_{X}\right) $ on $X;$ which
implies that the class $\mathrm{cl}\left( e_{1}\left[ p,p\right] \right) $
is a section of this sheaf. We denote this section by $\mathrm{ob}^{0}\left(
p\right) \in \Gamma \left( X,\mathcal{T}^{2}\left( \mathcal{O}_{X}\right)
\right) ;$ it can be called local obstruction to extension of $p$.

(ii) Suppose that the class $\mathrm{ob}^{0}\left( p\right) $ vanishes. We
can choose a sufficiently fine open covering $\left\{ Z_{\alpha },\alpha \in 
\mathcal{A}\right\} $ of $X$ and open sets $Y_{\alpha },X_{\alpha }$ for any 
$\,\alpha \in \mathcal{A}$ such that $Z_{\alpha }\Subset Y_{\alpha }\Subset
X_{\alpha }$ and for any index $\alpha $ a symmetric bounded operator $%
p_{2\alpha }:\Gamma \left( Y_{\alpha },\mathcal{O}_{X}\right) \otimes \Gamma
\left( Y_{\alpha },\mathcal{O}_{X}\right) \rightarrow \Gamma \left(
Z_{\alpha },\mathcal{O}_{X}\right) ,$ that fulfils (\ref{5}) for $a,b\in
\Gamma \left( Y_{\alpha },\mathcal{O}_{X}\right) $. We may assume also that
the restriction morphism $\Gamma \left( Y_{\alpha },\mathcal{O}_{X}\right)
\rightarrow \Gamma \left( Z_{\alpha },\mathcal{O}_{X}\right) $ is injective
and $Y_{\beta }\subset X_{\alpha }$ for all $\beta $ such that $Y_{\alpha
\beta }\doteq Y_{\alpha }\cap Y_{\beta }\neq \emptyset .$ Consider the
operator 
\begin{equation*}
q_{\alpha \beta }=p_{2\alpha }-p_{2\beta }:\Gamma \left( X_{\alpha },%
\mathcal{O}_{X}\right) \otimes \Gamma \left( X_{\alpha },\mathcal{O}%
_{X}\right) \rightarrow \Gamma \left( Z_{\alpha },\mathcal{O}_{X}\right) .
\end{equation*}
We have $\partial q_{\alpha \beta }=\partial p_{2\alpha }-\partial p_{2\beta
}=\left[ p,p\right] -\left[ p,p\right] =0.$ Following arguments of Lemma \ref
{d} this yields that $q_{\alpha \beta }$ is skew-symmetric derivation and
therefore can be uniquely extended to a derivation of sheaves $q_{\alpha
\beta }:\mathcal{O}_{X}\otimes \mathcal{O}_{X}\rightarrow \mathcal{O}_{X}$
defined in $Y_{\alpha \beta },$ that is $q_{\alpha \beta }\in \Gamma \left(
Y_{\alpha \beta },\mathcal{Q}_{X}\right) .$ Uniqueness implies that $%
q_{\beta \alpha }=-q_{\alpha \beta }.$

(iii) Consider the \v{C}ech 1-cochain $q\doteq \left\{ q_{\alpha \beta
}\right\} $ on the covering $\left\{ Y_{\alpha }\right\} .$ This a cochain
with values in the sheaf $\mathcal{Q}_{X}$ and obviously $\delta q=0.$%
Therefore the cohomology class $\mathrm{cl}\left( q\right) \in H^{1}\left( X,%
\mathcal{Q}_{X}\right) $ is well defined. If we choose some mappings $%
p_{2\alpha }^{\prime }$ instead of $p_{2\alpha },$we get the same class $%
\mathrm{cl}\left( q\right) .$ If $\mathrm{cl}\left( q\right) =0,$we have $%
q_{\alpha \beta }=r_{\alpha }-r_{\beta }.$ Replacing $p_{2\alpha }$ by $%
p_{2\alpha }^{\prime }=p_{2\alpha }-r_{\alpha }$ yields the equation $%
p_{2\alpha }^{\prime }=p_{2\beta }^{\prime }$ in $X_{\alpha \beta }.$ This
implies that three exists a global mapping $p^{\prime }:\mathcal{O}%
_{X}\times \mathcal{O}_{X}\rightarrow \mathcal{O}_{X}$ such that $p^{\prime
}|X_{\alpha }=p_{\alpha }^{\prime }.$ The class $\mathrm{ob}^{1}\left(
p\right) =\mathrm{cl}\left( q\right) $ can be called the global first
obstruction to extemsiom of $p$.

\begin{conclusion}
For an arbitrary complex analytic space $X$ there are defiend the
quadratic mapping $\mathrm{ob}^{0}:\Gamma \left( X,\mathcal{Q}_{X}\right)
\rightarrow \Gamma \left( X,\mathcal{T}^{2}\left( \mathcal{O}_{X}\right)
\right) ,$ $p\mapsto \mathrm{cl}\left( e_{2}\left[ p,p\right] \right) $ and
the homogeneous mapping $\mathrm{ob}^{1}:\mathrm{Ker\,ob}^{0}\rightarrow
H^{1}\left( X,\mathcal{Q}\left( X\right) \right) ,p\mapsto \mathrm{cl}\left(
q\right) .$ A Poisson bracket $p$ admits an extension to a star-product $%
\mathrm{mod}\,\left( \lambda ^{3}\right) $, if and only if $\mathrm{ob}%
^{0}p=0$ and $\mathrm{ob}^{1}p=0$.
\end{conclusion}

Further obstructions can be analyzed on the same lines. The obstruction
theory for deformations of complex analytic spaces looks very alike with the sheaf $\mathcal{%
Q}_{X}$ replaced by the sheaf $\mathcal{T}^{1}\left( \mathcal{O}_{X}\right) $
of Tyurina cohomology.

\section{Poisson brackets on K3-surfaces}

According to Kodaira's classification, a K3-surface is a compact analytic
2-manifold $X$ with the trivial canonical bundle $\Omega ^{2}\left( X\right) 
$. Note that for any surface the sheaf $\mathcal{Q}_{X}=\wedge ^{2}\mathcal{T%
}_{X}$ is dual to $\Omega ^{2}\left( X\right) $ and for any K3-surface
this sheaf is trivial which implies that $\mathrm{dim}\Gamma \left( X,%
\mathcal{Q}_{X}\right) =1.$ This means that there exists only one nontrivial
global skew-symmetric mapping $p:\mathcal{O}_{X}\otimes \mathcal{O}%
_{X}\rightarrow \mathcal{O}_{X}.$ Construct this mapping and show that this
is a global Poisson bracket. A simplest realization of a projective K3%
-surface is a quartic surface $X_{4}=\left\{ f_{4}=0\right\} $ in $\mathbb{CP%
}^{3}.$ Choose homogeneous coordinates $z_{0},z_{1},z_{2},z_{3}$ in $\mathbb{%
CP}^{3}$ and set 
\begin{equation*}
p_{ijk}\left( a,b\right) =\varepsilon _{ijkl}z_{l}^{-1}\det \left( 
\begin{array}{ccc}
\partial _{i}f & \partial _{j}f & \partial _{k}f \\ 
\partial _{i}a & \partial _{j}a & \partial _{k}a \\ 
\partial _{i}b & \partial _{j}b & \partial _{k}b
\end{array}
\right)
\end{equation*}
where $\partial _{i}=\partial /\partial z_{i}$ and $\left( i,j,k,l\right) $
is a permutation of $\left( 0,1,2,3\right) ,\varepsilon _{ijkl}$ is the sign
of this permutation. This mapping bracket is well defined in the domain $%
z_{l}\neq 0.$

\begin{proposition}
The mappings $J_{ijk}$ agree on the common domains.
\end{proposition}

\textsc{Proof.} Express the first column in the determinant by means of the
equations 
\begin{equation*}
\sum z_{j}\partial _{j}a=\sum z_{j}\partial _{j}b=0,\,\sum z_{j}\partial
_{j}f=4f
\end{equation*}
We obtain $p_{ijk}\left( a,b\right) =-\varepsilon
_{ijkl}z_{k}^{-1}p_{ijl}\left( a,b\right) .$ Taking in account that $%
-\varepsilon _{ijkl}=\varepsilon _{ijlk},$ we find that $p_{ijk}=p_{ijl}$. $%
\blacktriangleright $

There are just two more realizations of K3-surfaces as projective complete
intersections 
\begin{equation*}
X_{32}=\left\{ f_{3}=g_{2}=0\right\} \subset \mathbb{CP}^{4},\;X_{222}=\left%
\{ f_{2}=g_{2}=h_{2}=0\right\} \subset \mathbb{CP}^{5}
\end{equation*}
For the surface $X_{32}$ we define the bracket by means of homogeneous
coordinates $z_{0},...,z_{4}:$ 
\begin{equation*}
p_{ijkl}\left( a,b\right) =\varepsilon z_{m}^{-1}\frac{\partial \left(
f,g,a,b\right) }{\partial \left( z_{i},z_{j},z_{k},z_{l}\right) },
\end{equation*}
where $\left( i,j,k,l,m\right) $ is a permutation of $\left(
0,1,2,3,4\right) $\ and $\varepsilon $ is the sign of the permutation. For
the surface $X_{222}$ we take the similar bilinear form $p_{K}\left(
a,b\right) =J_{K}\left( f,g,h,a,b\right) $. By the same arguments we check
that the global bilinear mapping $p_{X}:\mathcal{O}_{X}\times \mathcal{O}%
_{X}\rightarrow \mathcal{O}_{X}$ is well defined in both cases.

\begin{corollary}
The mapping $p$ is a global Poisson bracket for any surface $%
X=X_{4},X_{32},X_{222}$ as above. The first obstruction vanishes for this
bracket.
\end{corollary}

The first fact follows from Proposition \ref{pb}. The second one is a
corollary of the equations $\mathcal{T}^{2}\left( X\right) =0$ and $%
H^{1}\left( X,\mathcal{O}_{X}\right) =0,$ since of arguments of Sec. 5. Note
that any K3-surface which is a projective complete intersection belongs to
one of the above types, whereas these types contain also singular algebraic
surfaces for special choice of polynomials $f,g,h.$ The above formulae are
well defined anyway.

\section{Quantization and deformation together}

The above examples have a special feature. First remind that the Poisson
bracket is nothing else as an infinitesimal quantization of the family of $%
K3 $-surfaces $X$ with the base $D_{2},$ whose support is one point $\cdot $
and the structure algebra is $\mathcal{O}\left( D_{2}\right) =R_{1}/\left(
\lambda ^{3}\right) .$ This means just that the mapping $p_{2}$ can be
found. In the same time, the family $\left\{ X\right\} $ itself is a flat
deformation in the category of compact complex spaces for all three types $%
X=X_{4},\,X_{32},\,X_{222}$ of algebraic surfaces. The deformation
parameters are just the coefficients of the polynomials $f,g,h$
respectively. The coefficients of $f,g,h$ run over the corresponding
projective spaces. For surfaces of types $X_{32}$ and $X_{222}$ the
coefficients must avoid some Zariski closed subspaces, where the dimension
of the common set of zero jumps up. The total space of coefficients is
redundant: one can take any subspace $S$ that is transversal to orbits of
the group of projective transformation acting on the polynomials. The
minimal local subspace $S$ is an open subset in $\,\mathbb{C}^{19}$. The
join of the deformation and of the quantization yields an associative
deformation of $X$ with base $S\times D_{2}$.

Note that the above family $\{X\}$ is versal in neither of the cases, since $%
\mathrm{\dim \,}T^{1}\left( X\right) =20$ and the base of a minimal versal
(commutative) deformation of $X$ is of dimension $20$. Deformation of an
algebraic surface $X$ in one special direction is still K3, but is no more
algebraic surface. The Poisson bracket $p$ as above must have an extension
to this versal deformation, but no explicit formula is known.

\section{Calculation of (co)homology of analytic algebras}

We prove here the following statements:

\begin{theorem}
\label{coh}For an arbitrary analytic algebra $A$ and any $n\geq 0$ the analytic (co)homology 
$A$-module $\mathrm{Hoch}_{n}\left( A,A\right) $ and $\mathrm{%
Hoch}^{n}\left( A,A\right) $ has finite set of generators.%
\newline
Moreover, for any complex analytic space $\left( X,\mathcal{O}_{X}\right) $
the sheaves $\mathrm{Hoch}_{n}\left( \mathcal{O}_{X}\mathcal{,O}_{X}\right) $
and $\mathrm{Hoch}^{n}\left( \mathcal{O}_{X}\mathcal{,O}_{X}\right) $ of
analytic (co)homology are coherent for $n\geq 0.$
\end{theorem}

\textsc{Proof.} The case $n=0$ is trivial since the (co)homology is
isomorphic to $A;$ we assume now that $n>0.$ The module of differentials on
an analytic algebra $A$ is defined to be $\Omega \left( A\right) =\Delta
/\Delta ^{2},$ where $\Delta $ is the ideal in the algebra $\mathrm{Ca}%
_{1}\left( A\right) =A^{\odot 2}$ that consists of elements $a$\ such that $%
\partial _{0}c\left( z\right) \doteq c\left( z,z\right) =0$ ; $\Delta ^{2}$
is the square of this ideal. For an arbitrary $a\in \mathrm{Ca}_{1}\left(
A\right) $ the image of $a\left( z,w\right) -a\left( z,z\right) $ in $\Omega
\left( A\right) $ is denoted by $\mathrm{d}a.$ In algebraic terms, we can
write $\mathrm{d}a=a-\sigma \partial _{0}a$ $\left( \mathrm{mod}\,\Delta
^{2}\right) ,$ where $\sigma b=b\left[ 1\right] ,b\in \mathrm{Ca}_{0}\left(
A\right) .$ The ideal $\Delta ^{2}$ is contained in the image of $\partial ,$
since $ab=\partial c$ for any $a,b\in \Delta $ and $c\left(
z^{0},z^{1},z^{2}\right) =-a\left( z^{0},z^{1}\right) b\left(
z^{0},z^{2}\right) .$ Therefore there is a morphism $\Omega \left( A\right)
\rightarrow \mathrm{Hoch}_{1}\left( A,A\right) .$

A free graded commutative analytic algebra $B=\oplus _{k\geq 0}B^{k}$ is, by
definition, a free graded commutative extension $A\left[ e_{1},e_{2},...%
\right] $ of an analytic algebra $A$ by means of homogeneous elements $%
e_{1},e_{2},...$ of strictly positive degrees such that for any $k$ the
number of elements of degree $k$ is finite. Each term $B^{k}$ is $A$-module
of finite type. The above definitions are generalized for the graded case;
the chain and cochain complexes acquire addition grading.

We formulate a version of the Hochschild-Kostant-Rosenberg\ theorem \cite
{HKR} for graded analytic algebras. For an arbitrary $n,$ there is the
canonical morphism of $B$-modules $\omega _{n}:\wedge ^{n}\Omega \left(
B\right) \rightarrow \mathrm{Ca}_{n}\left( B\right) $ defined for $%
b_{1},...,b_{n}\in \Delta \left( B\right) $ by 
\begin{equation}
\mathrm{d}b_{1}\wedge ...\wedge \mathrm{d}b_{n}\rightarrow \sum_{\pi \in
S_{n}}\left( -1\right) ^{\varepsilon \left( \pi \right) }\left[ b_{\pi
\left( 1\right) }|...|b_{\pi \left( n\right) }\right] .  \label{15}
\end{equation}
The sum is taken over all permutations $\pi $ and $\varepsilon \left( \pi
\right) =\sum \left| b_{\pi \left( i\right) }\right| \left| b_{\pi \left(
j\right) }\right| +1,$ where the sum is taken over all pairs $i<j$ such that 
$\pi \left( i\right) >\pi \left( j\right) $ and $\left| b\right| $ means the
degree of $b.$ We read the bracket as follows 
\begin{equation*}
\left[ b_{1}|...|b_{n}\right] \left( z^{0},...,z^{n}\right) =b_{1}\left(
z^{0},z^{1}\right) ...b_{n}\left( z^{0},z^{n}\right) 
\end{equation*}
when we interpret elements $b_{1},...b_{n}$ as functions of coordinates $%
z^{j},j=1,...,n$ evaluated in $\mathbb{C}^{n}$ by means of a surjection $%
R_{n}\rightarrow A$ and as polynomials of homogeneous generators $e$ of the
algebra $B.$ We have 
\begin{equation*}
\partial \left[ b_{1}|...|b_{n}\right] =\sum_{i=0}^{n-1}\left( -1\right) ^{i}%
\left[ b_{1}|...\left| b_{i}b_{i+1}\right| |b_{n}\right] 
\end{equation*}
and each term appear twice in the right-hand side with opposite signs which
implies that the sum belongs to $\mathrm{Ker\,}\partial $. It is contained
in $\mathrm{Im}\,\partial ,$ if $b_{j}\in \Delta ^{2}\left( B\right) $ for,
at least, one $j.$ Therefore (\ref{15}) generates the morphism of $B$%
-modules 
\begin{equation}
\omega _{n}:\wedge _{B}^{n}\Omega \left( B\right) \rightarrow \mathrm{Hoch}%
_{n}\left( B,B\right)   \label{9}
\end{equation}
The mapping (\ref{9}) induces by duality the morphism 
\begin{equation}
\omega ^{n}:\mathrm{Hoch}^{n}\left( B,B\right) \rightarrow \wedge _{B}^{n}%
\mathrm{Der}\left( B,B\right) ,  \label{16}
\end{equation}
where $\mathrm{Der}\left( B,B\right) =\mathrm{Hom}_{B}\left( \Omega \left(
B\right) ,B\right) $ is the module of derivations $t:B\rightarrow B.$

\begin{theorem}
\label{HKRM}Let $B$ be the free graded commutative algebra over a regular
analytic algebra $R_{n}$. Then the morphisms (\ref{9}) and (\ref{16}) are
bijections for any $n\geq 1$, where $\mathrm{Hoch}_{\ast }\left( B,B\right)
,\;\mathrm{Hoch}^{\ast }\left( B,B\right) $ mean the analytic Hochschild
(co)homology as above.
\end{theorem}

We skip a proof, which is inspired by arguments of the original paper \cite
{HKR}.

\textsc{Proof} of Theorem \ref{coh}. Let $A=R/I$ be an analytic algebra and $%
\left( R^{\ast },\rho \right) $ be its Tyurina resolvent. This means that $%
R^{\ast }=\oplus _{k\geq 0}R^{k}$ is a free graded commutative analytic
algebra; The differential $\rho $ is a derivation in $R^{\ast }$ of degree $%
-1,$ that is $\rho \left( ab\right) =\rho \left( a\right) b+\left( -1\right)
^{\left| a\right| }a\rho \left( b\right) .$ Any analytic algebra has a
Tyurina resolvent, see \cite{Pa} (where the grading in $R$ is taken
negative). Take the analytic chain complex $\mathrm{Ca}_{\ast }\left(
R^{\ast }\right) $ of the graded analytic algebra $R^{\ast }$ and define the
differential $\rho $ which is defined on products of homogeneous elements as
follows 
\begin{eqnarray*}
\rho \left( a_{0}\left[ a_{1}|a_{2}|...|a_{n}\right] \right)  &=&\rho \left(
a_{0}\right) \left[ a_{1}|a_{2}|...|a_{n}\right] +\left( -1\right) ^{\left|
a_{0}\right| }a_{0}\left[ \rho \left( a_{1}\right) |a_{2}|...|a_{n}\right] 
\\
&&+\left( -1\right) ^{\left| a_{0}\right| +\left| a_{1}\right| }a_{0}\left[
a_{1}|\rho \left( a_{2}\right) |...|a_{n}\right]  \\
&&+...+\left( -1\right) ^{\left| a_{0}\right| +\left| a_{1}\right|
+...+\left| a_{n-1}\right| }a_{0}\left[ a_{1}|a_{2}|...|\rho \left(
a_{n}\right) \right] ,
\end{eqnarray*}
where $a_{0},...,a_{n}\in R^{\ast }.$ This operator $\rho $ has extension to
the analytic tensor product $\mathrm{Ca}_{n}\left( R^{\ast }\right) =R^{\ast
\odot \left( n+1\right) };$ this extension is a bounded differential for any 
$n\geq 1.$ The grading $\left| a_{0}\left[ a_{1}|a_{2}|...|a_{n}\right]
\right| =\sum \left| a_{j}\right| $ in the tensor product induces a grading
in $\mathrm{Ca}_{n}\left( R^{\ast }\right) .$

\begin{lemma}
\label{s}The complex $\left( \mathrm{Ca}_{n}\left( R^{\ast }\right) ,\rho
\right) $ is exact and splits in positive degrees, and we have 
\begin{equation*}
H_{0}\left( \mathrm{Ca}_{n}\left( R^{\ast }\right) ,\rho \right) \cong 
\mathrm{Ca}_{n}\left( A\right) 
\end{equation*}
\end{lemma}

See \cite{P} for a proof. Introduce the total grading $a_{0}\left| \left[
a_{1}|a_{2}|...|a_{n}\right] \right| =\left| a_{1}\right| +...+\left|
a_{n}\right| +n$ in this complex $\mathrm{Ca}_{\ast }\left( R^{\ast }\right) 
$. The chain differential $\partial $ has degree $-1$ as well as $\rho .$
These differentials commute, since $\rho $ is a derivation. We have a
bicomplex $\left( \mathrm{Ca}_{\ast }\left( R^{\ast }\right) ,\partial ,\rho
\right) $, where both grading are positive. There are two spectral sequences 
$E,\,E^{\prime }$ that converge to the homology of the total complex $%
\mathrm{Tot}\left( \mathrm{Ca}_{\ast }\left( R^{\ast }\right) \right) $.

The sequence $E$ is generated by the filtration $F^{n}=\oplus _{m\leq n}%
\mathrm{Ca}_{m}\left( R^{\ast }\right) ,n=0,1,2,....$ We have $E_{0}=\mathrm{%
Ca}_{\ast }\left( R^{\ast }\right) $, $d_{0}=\rho $; by Lemma \ref{s} $%
E_{1}^{\ast 0}\cong \mathrm{Ca}_{\ast }\left( A\right) $ and $\,E_{1}^{\ast
k}=0$ for $k>0.$ This yields $E_{\infty }=E_{2}\cong H_{\ast }\left( \mathrm{%
Ca}_{\ast }\left( A\right) ,\partial \right) =\mathrm{Hoch}_{\ast }\left(
A,A\right) .$ Take another filtration $F^{\prime k}=\oplus _{i\leq k}\mathrm{%
Ca}_{\ast }\left( R^{\ast }\right) ^{i},$ where $\mathrm{Ca}_{\ast }\left(
R^{\ast }\right) ^{i}$ means the subspace of elements $a,\left| a\right| =i.$
Let $E^{\prime }$ be the corresponding spectral sequence. We have again $%
E_{1}^{\prime }=\mathrm{Ca}_{\ast }\left( R^{\ast }\right) ,$ but $%
d_{1}=\partial .$ It follows that $E_{2}^{\prime }=H_{\ast }\left( \mathrm{Ca%
}_{\ast }\left( R^{\ast }\right) ,\partial \right) .$ By Theorem \ref{HKRM}
we have 
\begin{equation*}
E_{2}^{\prime n}=H_{n}\left( \mathrm{Ca}_{\ast }\left( R^{\ast }\right)
,\partial \right) \cong \wedge ^{n}\Omega \left( R^{\ast }\right)
,\,n=1,2,...
\end{equation*}
The differential $d_{2}$ in $E_{2}^{\prime }$ is generated by the
differential $\rho$: 
\begin{equation*}
d_{2}\left( a_{0}\mathrm{d}a_{1}\wedge ...\wedge \mathrm{d}a_{n}\right)
=\rho \left( a_{0}\right) \mathrm{d}a_{1}\wedge ...\wedge \mathrm{d}a_{n}
\end{equation*}
\begin{equation*}
+\left( -1\right) ^{\left| a_{0}\right| }a_{0}\mathrm{d}\left( \rho
a_{1}\right) \wedge \mathrm{d}a_{2}\wedge ...\wedge \mathrm{d}a_{n}
+...+\left( -1\right) ^{\left| a_{0}\right| +...+\left| a_{n-1}\right| }a_{0}%
\mathrm{d}a_{1}\wedge ...\wedge \mathrm{d}\rho \left( a_{n}\right) . 
\end{equation*}
The module $\wedge ^{n}\Omega \left( R^{\ast }\right) $ is finitely
generated in each grading and the differential $d_{2}$ is a morphism of $A$%
-modules, as well as $d_{3},...$ Therefore the term $E_{\infty }^{\prime }$
is finitely generated in each degree. The equation $E_{\infty }^{\prime
}\cong H_{\ast }\left( \mathrm{Ca}_{\ast }\left( A\right) ,\partial \right) $
implies the statement of the theorem for the analytic Hochschild homology $%
\mathrm{Hoch}_{\ast }\left( A,A\right) $.

To calculate the analytic cohomology we consider the bicomplex $\mathrm{Ca}%
^{\ast }\left( R^{\ast }\right) ,$ where the first differential $\partial
^{\ast }$ is the standard one and the second differential $\rho ^{\ast }$is
defined by 
\begin{equation*}
\rho ^{\ast }h\left( a_{1}|...|a_{n}\right) =\rho \left( h\left(
a_{1}|...|a_{n}\right) \right) -\left( -1\right) ^{\left| h\right| }h\left(
\rho \left( a_{1}|...|a_{n}\right) \right) ,
\end{equation*}
where $\left| h\right| =\left| h\left( a\right) \right| -\left| a\right| $
means that degree of an operator $h:\mathrm{Ca}_{\ast }\left( R^{\ast
}\right) \rightarrow R^{\ast }.$ The differentials commute and there are two
spectral sequences that converges to the same limit. The first one is
induced by the grading $F^{n}=$ $\oplus _{m\geq n}\mathrm{Ca}^{m}\left(
R^{\ast }\right) .$ We have $E_{0}\cong \mathrm{Ca}^{\ast }\left( R^{\ast
}\right) $ and the differential $d_{0}$ is generated by $\rho ^{\ast }.$ It
follows from Lemma \ref{s} that $E_{1}^{\ast 0}\cong \mathrm{Ca}^{\ast
}\left( A\right) ,\,d_{1}=\partial ^{\ast }$ and $E_{1}^{\ast k}=0,\,k\neq
0. $ Therefore $E_{2}\cong \mathrm{Hoch}^{\ast }\left( A,A\right) .$ The
second filtration is formed by the subcomplexes 
\begin{equation*}
F^{\prime k}=\oplus _{i\leq k}\mathrm{Ca}^{\ast }\left( R^{\ast }\right)
^{i},\,k\in\mathbb{Z},
\end{equation*}
where the superscript means the degree as above. We have $E_{1}^{\prime
}\cong \mathrm{Ca}^{\ast }\left( R^{\ast }\right) $ and $d_{2}=\partial
^{\ast }$ The term $\,E_{2}^{\prime }$ is isomorphic to the analytic
cohomology $\mathrm{Hoch}^{\ast }\left( R^{\ast },R^{\ast }\right) .$
According to Theorem \ref{HKRM}, we have an isomorphism of $R^{\ast }$%
-modules $E_{2}^{\prime \ast }\cong \wedge ^{\ast }\mathrm{Der}\left(
R^{\ast },R^{\ast }\right) $ and $d_{2}$ is a differential $\wedge ^{\ast }%
\mathrm{Der}\left( R^{\ast },R^{\ast }\right) $ generated by the morphism $%
\mathrm{d}\rho ^{\ast }$ which is defined by 
\begin{eqnarray}
\mathrm{d}\rho ^{\ast }\left( t_{1}\wedge ...\wedge t_{n}\right) \left[
a_{1}|...|a_{n}\right] &=&\rho \left( t_{1}\left( a_{1}\right)
...t_{n}\left( a_{n}\right) \right)  \label{13} \\
&&-\left( -1\right) ^{\left| t\right| }\left( t_{1}\wedge ...\wedge
t_{n}\right) \left( \rho \left[ a_{1}|...|a_{n}\right] \right) ,  \notag
\end{eqnarray}
where $\left| t\right| \doteq \left| t_{1}\right| +...+\left| t_{n}\right| .$
By means of Lemma \ref{s} we can define a homotopy of the complex $\left(
\wedge ^{\ast }\mathrm{Der}\left( R^{\ast },R^{\ast }\right) ,d_{2}\right) $
to the complex $\left( \wedge ^{\ast }\mathrm{Der}\left( R^{\ast },A\right) ,%
\mathrm{d}\rho ^{\ast }\right) ,$ where the differential $\mathrm{d}\rho
^{\ast }$ looks like (\ref{13}), but the first term in the right-hand side
vanishes. This complex is $A$-module of finite type in each degree.
Therefore $E_{3}^{\prime }\cong H_{\ast }\left( \wedge ^{\ast }\mathrm{Der}%
\left( R^{\ast },A\right) ,\mathrm{d}\rho ^{\ast }\right) $ is also a $A$%
-module of finite type in each degree. The same true for $E_{\infty
}^{\prime },$ which completes the proof. $\blacktriangleright $

\textbf{Remark 1. }The (co)homology can be calculated by means of the
respective spectral sequences $\left( E_{2}^{\prime },d_{2}\right) .$ In
particular, the Harrison cohomology of $A$ is isomorphic to the cohomology
of the subcomplex of $\mathrm{Ca}^{\ast }\left( A\right) $ of cochains
vanishing on all shuffle-products. It is isomorphic to the Tyurina
cohomology which yields the injective mapping $T^{\ast }\left( A\right)
\rightarrow \mathrm{Hoch}^{\ast +1}\left( A,A\right) ,$ (\cite{P}). It can
be shown that the image of this mapping coincides with the contribution of $%
E_{2}^{0\ast }\cong H^{\ast }\left( \mathrm{Ca}^{1}\left( A\right) ,\rho
\right) $ in the cohomology of the total complex $\mathrm{Ca}^{\ast }\left(
R^{\ast }\right) .$

\textbf{Remark 2. }Theorem \ref{coh} can be generalized for analytic
homology and cohomology of an analytic algebra $A$ with values in an
arbitrary $A$-module $M$ of finite type.

\textbf{Remark 3. }The similar statements hold for the algebraic Hochschild
\linebreak (co)homology of any polynomial algebra $A=\mathbf{k}\left[
x_{1},...,x_{n}\right] /I$ and for any formal algebra $\mathbf{k}\left[ %
\left[ x_{1},...,x_{n}\right] \right] /I$. Calculations for complete
intersection polynomial algebras are done by Fr\o nsdal and Kontsevich, \cite
{FK}.

\bigskip

\end{document}